# EFFICIENT RANDOMIZED-ADAPTIVE DESIGNS

By Feifang Hu,[1] Li-Xin Zhang[2] and Xuming He[3]

*University of Virginia, Zhejiang University and University of Illinois at Urbana–Champaign*

Response-adaptive randomization has recently attracted a lot of attention in the literature. In this paper, we propose a new and simple family of response-adaptive randomization procedures that attain the Cramer–Rao lower bounds on the allocation variances for any allocation proportions, including optimal allocation proportions. The allocation probability functions of proposed procedures are discontinuous. The existing large sample theory for adaptive designs relies on Taylor expansions of the allocation probability functions, which do not apply to nondifferentiable cases. In the present paper, we study stopping times of stochastic processes to establish the asymptotic efficiency results. Furthermore, we demonstrate our proposal through examples, simulations and a discussion on the relationship with earlier works, including Efron's biased coin design.

**1. Introduction.** Scientific and clinical investigations often demand that a new intervention or treatment be compared against a control group. Randomization is the preferred way of assigning participants to the control and the treatment groups. In some clinical trials, response-adaptive randomization procedures are desirable for ethical and efficiency reasons [Hu and Rosenberger (2006)]. The three main components of response-adaptive randomization are allocation proportion, efficiency (power), and variability.

The development of response-adaptive randomization and optimal allocation proportion was briefly reviewed in the introduction of Hu and Zhang

Received April 2008; revised July 2008.
[1]Supported by Grant DMS-03-49048 from the National Science Foundation (USA).
[2]Supported by Grant 10771192 from the National Natural Science Foundation of China.
[3]Supported by the U.S. National Science Foundation Award DMS-06-04229, National Natural Science Foundation of China Grant 10828102 and a Changjiang Visiting Professorship at Northeast Normal University in China.

*AMS 2000 subject classifications.* Primary 60F15, 62G10; secondary 60F05, 60F10.
*Key words and phrases.* Response-adaptive designs, biased coin design, clinical trial, urn model, doubly adaptive biased coin design, power.

(2004). In this paper, we denote *optimal allocation proportions* as allocation proportions that are derived from certain optimality criteria (see examples in Section 2). The issue of efficiency or power was discussed by Hu and Rosenberger (2003), who showed that the efficiency is a decreasing function of the variability induced by the randomization procedure for any given allocation proportion. More recently, Hu, Rosenberger and Zhang (2006) showed that there was an asymptotic lower bound on the variability of response-adaptive designs. A response-adaptive design that attains this lower bound will be said to be first-order *efficient*. Most response-adaptive randomization procedures in the literature are not first-order efficient. The only exception that is known to us is the drop–the-loser rule [Ivanova (2003)], but its applications are limited to urn allocation proportion (not generally optimal) and binary responses [see Hu, Rosenberger and Zhang (2006) for details].

In this paper, we address the open problem of finding an efficient and randomized procedure that can adapt to any desired allocation proportion (including optimal allocation proportions). Specifically, we propose a new family of efficient randomized-adaptive designs (ERADE) that is easy to implement in practice for both discrete and continuous responses. Under some mild conditions, we obtain the asymptotic normality and strong consistency of both the allocation proportions and the estimators of population parameters. Our results provide a solid foundation for the efficient randomized-adaptive designs and the related statistical inferences based on such designs. Several advantages of the procedure are also demonstrated through examples.

In the literature, asymptotic properties of adaptive designs are studied under continuous and differentiable allocation probability function by using Taylor expansion. When the allocation probability function is discrete, the commonly used techniques do not work anymore. The allocation probabilities of the ERADEs are discrete functions. In this paper, we overcome the difficulties of discontinuity by introducing a stopping time of a martingale process. This may provide a novel direction to study the properties of adaptive designs with discrete probability functions. Efron's biased coin design [Efron (1971)] is a special case of ERADE for balancing two treatments. Rosenberger and Lachin [(2002), page 49] showed, by simulation, that Efron's biased coin design is more efficient than adaptive biased coin designs proposed by Wei (1978) and Smith (1984). The asymptotic results in this paper show the efficiency of Efron's biased coin design theoretically.

The paper is organized as follows. In Section 2, we propose a simple family of response-adaptive randomization procedures with some important asymptotic properties. They are shown to attain the Cramer–Rao lower bound on the allocation variabilities. The efficient procedures are derived for several commonly used optimal allocation proportions. In Section 3, we conduct a simulation study to examine the properties of the proposed designs. Section 4 demonstrates the potential value of the ERADE through a careful analysis of the ECMO trials. We conclude the paper with several remarks and potential further research directions in Section 5. Technical details are given in the Appendix.

**2. Efficient randomized-adaptive designs.** We first describe the framework for the randomized-adaptive designs.

*General framework.* We consider two-treatment clinical trials. Suppose that the patients come to the clinical trial sequentially and respond to treatments without delay. After the first $m$ patients are being assigned to treatments and the responses observed, the $(m+1)$th patient will be assigned to treatment 1 with probability $p_{m+1}$ and treatment 2 with probability $1 - p_{m+1}$. The probability $p_{m+1}$ may depend on both the treatments assigned to and the responses observed of the previous $m$ patients. Let $X_{m,k}$ be the result of the $m$th assignment (i.e., $X_{m,k} = 1$ if the $m$th patient is assigned to treatment $k$, and 0 otherwise). We assume that the patient responses with each treatment are i.i.d. with the probability distribution indexed by $\boldsymbol{\theta}_k \in R^d$ ($k = 1, 2$). The cases of $d = 1$ (for binary response) or $d = 2$ (for normally distributed response) are typical.

Let $N_{m,k} = \sum_{j=1}^{m} X_{j,k}$ be the number of patients assigned to treatment $k$ in the first $m$ patients. We assume that the desired allocation proportion of patients assigned to each treatment is a function of $\boldsymbol{\Theta} = (\boldsymbol{\theta}_1, \boldsymbol{\theta}_2)$ [see Rosenberger et al. (2001) for a related discussion]. More specifically, the goal of allocation is to have $N_{m,1}/m \to v = \rho(\boldsymbol{\Theta})$, where $\rho(\cdot) : \mathcal{R}^{d \times 2} \to (0, 1)$ is an allocation proportion.

We further assume that the parameter estimate $\widehat{\boldsymbol{\theta}}_{m,k}$ of $\boldsymbol{\theta}_k$ based on an $m$-patient study has the Bahadur-type representation

$$(2.1) \qquad \widehat{\boldsymbol{\theta}}_{m,k} = N_{m,k}^{-1} \sum_{j=1}^{m} X_{j,k} \boldsymbol{\xi}_{j,k} + o(N_{m,k}^{-1/2}) \qquad \text{as } m \to \infty$$

for some i.i.d. sequences of random variables $\{\boldsymbol{\xi}_{j,k}, j = 1, \ldots\}$, where $\boldsymbol{\xi}_{m,k}$ is the response or a function of the response of the $m$th patient on treatment $k$.

In most applications, the response distributions belong to an exponential family. Then, we take $\boldsymbol{\theta}_k = E[\boldsymbol{\xi}_{1,k}]$, where $\boldsymbol{\xi}_{1,k}$ are the natural sufficient statistics, and $\widehat{\boldsymbol{\theta}}_{m,k}$ is the average of the observed sufficient statistics. In practice, we may start with $\boldsymbol{\Theta}_0 = (\boldsymbol{\theta}_{0,1}, \boldsymbol{\theta}_{0,2})$ as an initial estimate of $\boldsymbol{\Theta}$ and use

$$(2.2) \qquad \widehat{\boldsymbol{\theta}}_{m,k} = \frac{\sum_{j=1}^{m} X_{j,k} \boldsymbol{\xi}_{j,k} + \boldsymbol{\theta}_{0,k}}{N_{m,k} + 1}$$

to ensure a well-defined estimator even when no patients are assigned to treatment $k$. The initial estimate $\Theta_0$ is a guessed value of $\Theta$ or an estimate of $\Theta$ from early trials.

*Efficient randomized-adaptive design (ERADE).* To start, we assign $m_0$ (usually $m_0 = 2$) subjects to each treatment by using restricted randomization. Assume that $m$ ($m \geq 2m_0$) subjects have been assigned in the trial, and their responses are observed. Let $\widehat{\Theta}_m$ be the estimator defined in (2.1) based on the $m$ observations, and $\widehat{\rho}_m = \rho(\widehat{\Theta}_m)$. Then, our proposed ERADE assigns the $(m+1)$th patient to treatment 1 with probability

$$(2.3) \qquad p_{m+1} = \begin{cases} \alpha \widehat{\rho}_m, & \text{if } N_{m,1}/m > \widehat{\rho}_m, \\ \widehat{\rho}_m, & \text{if } N_{m,1}/m = \widehat{\rho}_m, \\ 1 - \alpha(1 - \widehat{\rho}_m), & \text{if } N_{m,1}/m < \widehat{\rho}_m, \end{cases}$$

where $0 \leq \alpha < 1$ is a constant that reflects the degree of randomization. See Remarks 2.1 and 3.1 for details of choosing $\alpha$.

It is important to note that the allocation probability in (2.3) is a discrete function. In the literature, allocation probability functions are generally continuous and differentiable [Hu and Rosenberger (2006)], which allows for Taylor expansions in the asymptotic analysis. We have to use a different approach to establish the asymptotic properties in the Appendix.

EXAMPLE 1 (*Efron's biased Coin design*). To balance patients in two treatment groups (with $\rho = 1/2$ as the target), our proposed ERADE assigns the $(m+1)$th patient to treatment 1 with probability

$$(2.4) \qquad p_{m+1} = \begin{cases} \alpha/2, & \text{if } N_{m,1}/m > 1/2, \\ 1/2, & \text{if } N_{m,1}/m = 1/2, \\ 1 - \alpha/2, & \text{if } N_{m,1}/m < 1/2. \end{cases}$$

The special case of $\alpha = 2/3$ is Efron's biased coin design [Efron (1971)].

REMARK 2.1. Efron's biased coin design played a very important role in the randomization literature. This idea has been developed further by Pocock and Simon (1975), among others, for randomization procedures that depend on prognostic variables. However, it is difficult to study its asymptotic properties because the allocation probability is not a continuous function. To overcome the theoretical difficulty, Wei (1978) and Smith (1984) proposed adaptive biased coin designs and generalized biased coin designs by using continuous allocation probability. Burman (1996) introduced the expected $p$-value deficiency to evaluate the performance of a given design. Based on Burman's studies, it is reasonable to choose $\alpha$ between 0.4 and 0.7. Also, Baldi Antognini and Giovagnoli (2004) modified Efron's procedure by using more flexible forms of probability functions.

EXAMPLE 2 (*Binary response*). In the case of binary response, we denote $P_1 = 1 - Q_1$ and $P_2 = 1 - Q_2$, which are the success probabilities of treatments 1 and 2, respectively, and $\xi_{j,k}$ ($j = 1, \ldots, m$ and $k = 1, 2$) as the corresponding responses. If $X_{j,k}$ is the treatment assignment indicator of $j$th patient, as defined in the general framework, then the natural choices of parameter estimates are

$$\hat{P}_{m,k} = \frac{\sum_{j=1}^m X_{j,k}\xi_{j,k} + 0.5}{N_{m,k} + 1}, \qquad k = 1, 2.$$

To target $\rho(P_1, P_2)$ as the allocation proportion, let $\hat{\rho}_m = \rho(\hat{P}_{m,1}, \hat{P}_{m,2})$. Then, the ERADE assigns the $(m+1)$th patient to treatment 1 with probability $p_{m+1}$, as defined in (2.3). We now list three specific allocation proportions:

(i) The urn allocation proportion, given in Zelen (1969), Wei and Durham (1978) and Ivanova (2003),

$$(2.5) \qquad v_1 = \rho(P_1, P_2) = (1 - P_1)/(2 - P_1 - P_2);$$

(ii) The optimal allocation proportion discussed in Rosenberger et al. (2001),

$$(2.6) \qquad v_2 = \rho(P_1, P_2) = \sqrt{P_1}/(\sqrt{P_1} + \sqrt{P_2});$$

(iii) The Neyman allocation proportion discussed by Jennison and Turnbull (2000),

$$(2.7) \qquad v_3 = \rho(P_1, P_2) = \sqrt{P_1 Q_1}/(\sqrt{P_1 Q_1} + \sqrt{P_2 Q_2}).$$

EXAMPLE 3 (*Gaussian response*). We compare two treatments with responses $\xi_1 \sim N(\mu_1, \tau_1^2)$ and $\xi_2 \sim N(\mu_2, \tau_2^2)$, respectively. Let $\xi_{j,k}$, $X_{j,k}$ and $N_{m,k}$ be defined as in Example 2. Then, the maximum likelihood estimates are

$$\hat{\mu}_{m,k} = \frac{\sum_{j=1}^m X_{j,k}\xi_{j,k}}{N_{m,k}} \quad \text{and} \quad \hat{\tau}_{m,k}^2 = \frac{\sum_{j=1}^m X_{j,k}(\xi_{j,k} - \hat{\mu}_{m,k})^2}{N_{m,k}}, \qquad k = 1, 2.$$

One can show that the MLEs satisfy the Bahadur-type representation in (2.1). We can target any given allocation proportion $\rho(\mu_1, \mu_2, \tau_1, \tau_2)$ by using its estimator, $\hat{\rho}_m = \rho(\hat{\mu}_{m,1}, \hat{\mu}_{m,2}, \hat{\tau}_{m,1}, \hat{\tau}_{m,2})$, and the ERADE assigns the $(m+1)$th patient to treatment 1 with probability $p_{m+1}$ as defined in (2.3). Three specific examples are as follows:

(i) The optimal allocation proportion given by Zhang and Rosenberger (2006),

$$(2.8) \qquad v_4 = \rho(\mu_1, \mu_2, \tau_1, \tau_2) = \frac{\tau_1\sqrt{\mu_2}}{\tau_1\sqrt{\mu_2} + \tau_2\sqrt{\mu_1}},$$

when both $\mu_1$ and $\mu_2$ are nonnegative;

(ii) The Neyman allocation proportion used by Jennison and Turnbull (2000),

$$v_5 = \rho(\mu_1, \mu_2, \tau_1, \tau_2) = \frac{\tau_1}{\tau_1 + \tau_2}; \tag{2.9}$$

(iii) The $D_A$-optimal allocation proportion (Gwise, Hu and Hu (2008)),

$$v_6 = \rho(\mu_1, \mu_2, \tau_1, \tau_2) = \frac{\tau_1^{4/3}}{\tau_1^{4/3} + \tau_2^{4/3}}. \tag{2.10}$$

REMARK 2.2. In Example 2, $v_2$ and $v_3$ are optimal allocation proportions, where $v_2$ is based on minimizing the expected number of treatment failures for fixed variance of Wald test [Rosenberger et al. (2001)], and $v_3$ is based on minimizing the variance of the Wald test [Jennison and Turnbull (2000)]. Also, $v_1$ is not an optimal allocation proportion, but it is the limiting proportion of the play–the-winner rule [Zelen (1969)] and the randomized play–the-winner rule [Wei and Durham (1978)]. In Example 3, $v_4$, $v_5$ and $v_6$ are optimal allocation proportions from different optimal criteria.

By the results of Hu and Rosenberger (2003), the asymptotic properties, especially the asymptotic variance of the allocation, play an important role in the efficiency of the design. Now, we turn to the asymptotic properties of the proposed ERADE. First, we state two conditions.

CONDITION A. In the Bahadur-type representation (2.1), $\boldsymbol{\theta}_k = \mathsf{E}\boldsymbol{\xi}_{1,k}$ and $\mathsf{E}\|\boldsymbol{\xi}_{1,k}\|^2 < \infty$ ($k = 1, 2$).

CONDITION B. The proportion function $\mathbf{y} \to \rho(\mathbf{y}): \mathcal{R}^{d \times 2} \to (0, 1)$ is a continuous function and is twice differentiable at $\boldsymbol{\Theta}$ with $\rho(\boldsymbol{\Theta}) = v$.

These two conditions are mild, and they are satisfied by all the allocation functions discussed in Examples 1–3. We now state the main asymptotic properties of the ERADE with the notation $\mathbf{V}_k = \mathrm{Var}(\boldsymbol{\xi}_{1,k})$ ($k = 1, 2$), and

$$\mathbf{V} = \mathrm{diag}\left(\frac{1}{v}\mathbf{V}_1, \frac{1}{1-v}\mathbf{V}_2\right) \quad \text{and} \quad \sigma^2 = \left(\frac{\partial \rho}{\partial \mathbf{y}}\bigg|_{\boldsymbol{\Theta}}\right)' \mathbf{V} \frac{\partial \rho}{\partial \mathbf{y}}\bigg|_{\boldsymbol{\Theta}}.$$

THEOREM 2.1. *Under Conditions A and B, we have, as $n \to \infty$,*

$$\begin{aligned}
|N_{n,1} - n\widehat{\rho}_n| &= o_P(\sqrt{n}) \quad \text{and} \\
|N_{n,1} - n\widehat{\rho}_n| &= O(\sqrt{n \log \log n}) \qquad a.s.
\end{aligned} \tag{2.11}$$

*Further, we have*

(2.12) $$\sqrt{n}(\hat{\mathbf{\Theta}}_n - \mathbf{\Theta}) \xrightarrow{D} \mathrm{N}(0, \mathbf{V}),$$

(2.13) $$\sqrt{n}(N_{n,1}/n - v) \xrightarrow{D} \mathrm{N}(0, \sigma^2),$$

$$\hat{\mathbf{\Theta}}_n - \mathbf{\Theta} = O(\sqrt{n^{-1} \log \log n}) \qquad a.s.$$

*and*

(2.14) $$N_{n,1} - nv = O(\sqrt{n \log \log n}) \qquad a.s.$$

By Theorem 2.1 and the result of Hu, Rosenberger and Zhang (2006), we see that the proposed ERADE is indeed asymptotically efficient, so long as the parameter estimate $\hat{\mathbf{\Theta}}_n$ is efficient. Let $I_k$ be the Fisher information matrix for parameter $\boldsymbol{\theta}_k$ ($k = 1, 2$). Then, the efficiency result is formally summarized in the following theorem.

THEOREM 2.2. *Under Conditions A and B, if* $\mathrm{Var}(\boldsymbol{\xi}_{1,k}) = I_k^{-1}$, *then the asymptotic variance of* $N_{n,1}/\sqrt{n}$ *for the ERADE attains the Cramer–Rao lower bound*

(2.15) $$\left(\frac{\partial \rho}{\partial \mathbf{y}}\bigg|_{\mathbf{\Theta}}\right)' \mathrm{diag}((vI_1)^{-1}, ((1-v)I_2)^{-1}) \frac{\partial \rho}{\partial \mathbf{y}}\bigg|_{\mathbf{\Theta}}.$$

If the response distributions belong to the exponential family, Theorem 2.2 establishes the efficiency of the ERADE for any $\alpha \in [0, 1)$. It also explains why Efron's biased coin design works well for balancing two treatments ($\rho = 1/2$).

REMARK 2.3. Property (2.11) is a key component in the asymptotic analysis of the sample proportion $N_{n,1}/n$. It shows that the sample proportion can be approximated by an estimator $\hat{\rho}_n$ of the target proportion $\rho$. In the literature, where the allocation probability is generally a continuous and differentiable function of $\hat{\rho}_n$ and /or of the current sample proportion, the asymptotic properties of adaptive designs have been obtained by using Taylor expansion [cf., Wei (1978), Smith (1984) and Hu and Zhang (2004) among many others]. In those cases, the sample proportion cannot be well (efficiently) approximated by $\hat{\rho}_n$; therefore, the allocation variabilities do not attain the Cramer–Rao lower bound. When the allocation probability function is discrete, different techniques are needed and different asymptotic properties expected. We shall prove (2.11) by defining a stopping time of a martingale process to show that the difference between the sample proportion and the estimated target proportion in our model is small.

EXAMPLE 1 (*Continued*). Because $\rho = 1/2$ is a constant, the lower bound is 0. That is, $\sqrt{n}(N_{n,1}/n - 1/2) \to 0$ in probability. Therefore, Efron's biased coin design is a higher order approximation to balance two treatments. This explains the fact that Efron's biased coin design provides better balanced results than the adaptive biased coin design [Wei (1978)] and the generalized biased coin design [Smith (1984)].

EXAMPLE 2 (*Continued*). Both Conditions A and B are satisfied for binary response and the three allocations defined in (2.5), (2.6) and (2.7). By Theorem 2.1, we have, for the three allocations in (2.5), (2.6) and (2.7),

$$\sqrt{n}(N_{n,1}/n - v_i) \xrightarrow{D} N(0, \sigma_i^2), \qquad i = 1, 2, 3$$

with

$$\sigma_1^2 = \frac{Q_1 Q_2 (P_1 + P_2)^2}{(2 - P_1 - P_2)^3}, \tag{2.16}$$

$$\sigma_2^2 = \frac{Q_2 P_1^{1.5} + Q_1 P_2^{1.5}}{4\sqrt{P_1 P_2}(\sqrt{P_1} + \sqrt{P_2})^3} \tag{2.17}$$

and

$$\sigma_3^2 = \frac{[P_1 Q_1]^{1.5}(1 - 2P_2)^2 + [P_2 Q_2]^{1.5}(1 - 2P_1)^2}{4\sqrt{P_1 Q_1 P_2 Q_2}[\sqrt{P_1 Q_1} + \sqrt{P_2 Q_2}]^3}. \tag{2.18}$$

By Theorem 2.2, the variances in (2.16), (2.17) and (2.18) attain the Cramer–Rao lower bounds in their respective allocations.

EXAMPLE 3 (*Continued*). By using $\xi_k^* = (\xi_k, (\xi_k - \mu_k)^2)'$ as the responses, the conditions of Theorems 2.1 and 2.2 are satisfied. Thus, for the three allocations in (2.8), (2.9) and (2.10), we have $\sqrt{n}(N_{n,1}/n - v_i) \xrightarrow{D} N(0, \sigma_i^2)$ ($i = 4, 5, 6$) with

$$\sigma_4^2 = \frac{\tau_1 \tau_2 \sqrt{\mu_2} \sqrt{\mu_1}}{2(\tau_1 \sqrt{\mu_2} + \tau_2 \sqrt{\mu_1})^2}, \tag{2.19}$$

$$\sigma_5^2 = \frac{\tau_1 \tau_2}{2(\tau_1 + \tau_2)^2} \tag{2.20}$$

and

$$\sigma_6^2 = \frac{8(\tau_1 \tau_2)^{4/3}}{9(\tau_1^{4/3} + \tau_2^{4/3})^2}. \tag{2.21}$$

These variances in (2.19), (2.20) and (2.21) attain the Cramer–Rao lower bounds in their respective allocations.

TABLE 1
*Results for Case 1. Simulated and theoretical allocation proportions $(N_{n,1}/n)$ are given for different designs with their variances $\mathrm{nvar}(N_{n,1}/n)$ given in the parentheses. The target allocation proportion is $v_1 = (1 - P_2)/(2 - P_1 - P_2)$. The simulation used 1000 trials of $n = 100$*

| | | **ERADE with** | | | **DBCD** | |
| | **ERADE and DL** | $\alpha = 1/2$ | $\alpha = 2/3$ | **DL** | **with $\gamma = 2$** | |
| $P_1, P_2$ | Asymptotic | Simulated | Simulated | Simulated | Asymptotic | Simulated |
|---|---|---|---|---|---|---|
| 0.9, 0.7 | 0.75 (0.75) | 0.72 (0.69) | 0.72 (0.73) | 0.64 (0.39) | 0.75 (0.94) | 0.74 (0.91) |
| 0.9, 0.6 | 0.80 (0.48) | 0.78 (0.49) | 0.77 (0.51) | 0.69 (0.34) | 0.80 (0.61) | 0.78 (0.58) |
| 0.9, 0.5 | 0.83 (0.32) | 0.81 (0.34) | 0.80 (0.36) | 0.73 (0.24) | 0.83 (0.42) | 0.81 (0.39) |
| 0.9, 0.3 | 0.88 (0.16) | 0.85 (0.17) | 0.85 (0.18) | 0.79 (0.14) | 0.88 (0.22) | 0.86 (0.21) |
| 0.8, 0.8 | 0.50 (10.00) | 0.50 (0.75) | 0.50 (0.82) | 0.50 (0.51) | 0.50 (1.25) | 0.50 (1.16) |
| 0.8, 0.7 | 0.60 (0.72) | 0.59 (0.65) | 0.59 (0.63) | 0.57 (0.43) | 0.60 (0.91) | 0.60 (0.80) |
| 0.8, 0.6 | 0.67 (0.52) | 0.66 (0.52) | 0.65 (0.49) | 0.62 (0.35) | 0.67 (0.67) | 0.66 (0.63) |
| 0.7, 0.5 | 0.63 (0.35) | 0.62 (0.35) | 0.62 (0.36) | 0.60 (0.30) | 0.63 (0.47) | 0.62 (0.45) |
| 0.7, 0.3 | 0.70 (0.21) | 0.69 (0.20) | 0.69 (0.23) | 0.68 (0.18) | 0.70 (0.29) | 0.69 (0.28) |
| 0.6, 0.4 | 0.60 (0.24) | 0.60 (0.24) | 0.59 (0.25) | 0.59 (0.23) | 0.60 (0.34) | 0.60 (0.31) |
| 0.5, 0.5 | 0.50 (0.25) | 0.50 (0.22) | 0.50 (0.23) | 0.50 (0.21) | 0.50 (0.35) | 0.50 (0.33) |
| 0.5, 0.2 | 0.62 (0.13) | 0.61 (0.13) | 0.61 (0.15) | 0.61 (0.13) | 0.62 (0.20) | 0.61 (0.20) |
| 0.4, 0.3 | 0.54 (0.13) | 0.54 (0.13) | 0.54 (0.13) | 0.54 (0.13) | 0.54 (0.21) | 0.54 (0.20) |
| 0.2, 0.2 | 0.50 (0.06) | 0.50 (0.06) | 0.50 (0.06) | 0.50 (0.07) | 0.50 (0.13) | 0.50 (0.12) |

REMARK 2.4. In Theorem 2.1, the sample proportion $N_{n,1}/n$ converges to $v$ (target proportion) at the rate of $n^{-1/2}$. The asymptotic results in Theorems 2.1 and 2.2 do not depend on $\alpha$. This is because the allocation probability is a discrete function and the first order approximation does not depend on $\alpha$. In practice, we need to choose a suitable $\alpha$. We will discuss this in next section.

**3. Simulation study.** We conducted Monte Carlo simulations to examine the finite sample performance of the proposed ERADE. As shown by Hu and Rosenberger (2003), the efficiency (power) of a randomized design is a decreasing function of the variability of the design, so the variability of the allocation proportions is the main criterion in our investigation. We compared the doubly adaptive biased coin design (DBCD) of Hu and Zhang (2004), the drop–the-loser (DL) rule of Ivanova (2003), the randomized player–the-winner (RPW) rule of Wei and Durham (1978) and our proposed ERADE with $\alpha = 1/2$ or $2/3$. Here, we use $\alpha = 1/2$ or $2/3$, based on Burman's studies of Efron's biased coin design. Note that the DL rule and the RPW rule apply only to binary responses by targeting the urn allocation proportion, which is not optimal. On the other hand, the ERADE applies

TABLE 2
*Results for Case 2 with the target allocation proportion $v_2 = \sqrt{P_1}/(\sqrt{P_1} + \sqrt{P_2})$. Simulated and theoretical allocation proportions $(N_{n,1}/n)$ are given for different designs with their variances $\operatorname{nvar}(N_{n,1}/n)$ given in the parentheses. The simulation used 1000 trials of $n = 100$*

| | | **ERADE with** | | **DBCD** | |
| | **ERADE** | $\alpha = 1/2$ | $\alpha = 2/3$ | **with $\gamma = 2$** | |
| $P_1, P_2$ | Asymptotic | Simulated | Simulated | Asymptotic | Simulated |
|---|---|---|---|---|---|
| 0.9, 0.7 | 0.53 (0.02) | 0.53 (0.02) | 0.53 (0.03) | 0.53 (0.07) | 0.53 (0.07) |
| 0.9, 0.6 | 0.55 (0.03) | 0.55 (0.03) | 0.55 (0.04) | 0.55 (0.08) | 0.55 (0.08) |
| 0.9, 0.5 | 0.57 (0.04) | 0.57 (0.04) | 0.57 (0.06) | 0.57 (0.09) | 0.57 (0.11) |
| 0.9, 0.3 | 0.63 (0.09) | 0.64 (0.12) | 0.63 (0.13) | 0.63 (0.15) | 0.64 (0.18) |
| 0.8, 0.8 | 0.50 (0.02) | 0.50 (0.02) | 0.50 (0.03) | 0.50 (0.07) | 0.50 (0.07) |
| 0.8, 0.7 | 0.52 (0.02) | 0.52 (0.03) | 0.52 (0.03) | 0.52 (0.08) | 0.52 (0.08) |
| 0.8, 0.6 | 0.54 (0.03) | 0.54 (0.04) | 0.53 (0.04) | 0.54 (0.09) | 0.54 (0.09) |
| 0.7, 0.5 | 0.54 (0.05) | 0.54 (0.05) | 0.54 (0.06) | 0.54 (0.10) | 0.54 (0.12) |
| 0.7, 0.3 | 0.60 (0.09) | 0.60 (0.12) | 0.60 (0.13) | 0.60 (0.16) | 0.61 (0.21) |
| 0.6, 0.4 | 0.55 (0.07) | 0.55 (0.08) | 0.55 (0.09) | 0.55 (0.13) | 0.55 (0.15) |
| 0.5, 0.5 | 0.50 (0.06) | 0.50 (0.07) | 0.50 (0.08) | 0.50 (0.13) | 0.50 (0.14) |
| 0.5, 0.2 | 0.61 (0.17) | 0.62 (0.19) | 0.61 (0.22) | 0.61 (0.25) | 0.62 (0.30) |
| 0.4, 0.3 | 0.54 (0.12) | 0.53 (0.14) | 0.54 (0.16) | 0.54 (0.19) | 0.54 (0.25) |
| 0.2, 0.2 | 0.50 (0.25) | 0.50 (0.28) | 0.50 (0.33) | 0.50 (0.35) | 0.50 (0.48) |

to all types of responses and can be used to target any desired allocation proportion.

CASE 1. The data were generated from the binary responses based on Example 2. First, we consider the urn allocation proportion in (2.5), and the results are given in Table 1. By Hu and Rosenberger (2003), the RPW rule has much higher variability than the DBCD, so we do not include the results of the former here. For the ERADE, $m_0 = 2$ and $\alpha = 1/2$ or $\alpha = 2/3$ are considered. For the DL rule, we start with 5 type 1 balls, 5 type 2 balls and 1 type 0 ball. For the DBCD, we use the allocation function of Hu and Resenberger (2003), with $\gamma = 2$. The results in Table 1 are based on 1000 simulated trials with sample size 100. From Table 1, we see that the ERADEs with $\alpha = 1/2$ and $\alpha = 2/3$ performed better than other two procedures in general. The DBCD had larger variances than other procedures. Both the ERADE and the DBCD attained the target allocation proportion quite accurately, and the simulated variance agreed with their corresponding asymptotic variances quite well, except $P_1 = P_2 = 0.80$. When $P_1$ and $P_2$ were large and different, the DL rule did not converge to the target allocation proportion as fast as other procedures. However, its finite-sample variances were much smaller than the corresponding asymptotic variances. When both $P_1$ and $P_2$ were small, both the DL rule and the ERADE performed well.

CASE 2. When the target allocation is the optimal allocation proportion in (2.6), the DL rule no longer applies. The comparison between the ERADE and the DBCD are given in Table 2. Both designs converged well to the optimal allocation proportion, but the ERADE performed better with much smaller variability. Similar numerical results were obtained for the Neyman allocation proportion in (2.7), which are not reported here. We also considered the cases with sample sizes 50 and 200, but only to reach the same conclusions.

REMARK 3.1. It is easy to see that $\alpha$ is related to the randomness of the design. When $\alpha$ is smaller, the ERADE is more determined and could have smaller variability. To study this, we have run simulations (both Cases 1 and 2) for $\alpha = 1/8$, $1/4$, $1/2$, $2/3$ and $3/4$ with different sample sizes $n = 50$, 100 and 200. We found that the simulated results of $\alpha = 1/8$ and $1/4$ are very similar to the results of $\alpha = 1/2$ in terms of allocation proportion and its variability. The ERADE with $\alpha = 3/4$ has a slightly larger variability than others. Therefore, it is reasonable to recommend choosing $\alpha$ between 0.4 and 0.7. This agrees with Burman's (1996) study of biased coin designs.

REMARK 3.2. The simulated results in Tables 1 and 2 are based on binary responses from Example 2. We have also conducted some Monte Carlo simulations based on continuous responses from Example 3. We considered the doubly adaptive biased coin design (DBCD with $\gamma = 2$) of Hu and Zhang (2004) and our proposed ERADEs with $\alpha = 1/2$ and $2/3$. We obtained similar results as reported in Table 2 with sample size $n = 100$. All designs converged well to the target allocation proportion, but the ERADE performed better with smaller variability.

**4. Analysis of ECMO trials, where ERADE would make a difference.** Extracorporeal membrane oxygenation (ECMO) is an external system for oxygenating the blood based on techniques used in cardiopulmonary bypass technology developed for cardiac surgery. In the literature, there are three well-documented clinical trials on the evaluation of the clinical effectiveness of ECMO. These are the Michigan ECMO study [Bartlett et al. (1985)], the Boston ECMO study [Ware (1989)] and the UK ECMO trial [UK Collaborative ECMO Trials Group (1996)]. The Boston ECMO study is based on a two-stage design, so we restrict our discussion to the other two studies here.

In the Michigan ECMO study, the randomized player–the-winner rule was used to allocate patients to two treatments. Out of 12 patients, only one infant received conventional therapy. The first infant was assigned to ECMO and survived. The second was assigned to conventional therapy and died. Ten subsequent infants were randomized to ECMO and all survived. This trial provided very little information about survival rates in the same population treated with conventional therapy, because only one infant was assigned to the conventional therapy. It turned out that this patient was the most severely ill patient in the study. The high variability of the randomized player–the-winner rule was responsible for the failure of the Michigan ECMO study.

The UK ECMO trial used randomized allocation with equal proportions, and there were 93 patients in the ECMO and 92 in the conventional treatment for a total of 185. Prior to discharge from hospital, there were 28 deaths in the ECMO treatment and 54 deaths in the conventional treatment. We will use $P_1 = 65/93$ and $P_2 = 38/92$ as the estimated success probabilities of the ECMO and the conventional treatment, respectively. Taking those values as the success probabilities, we discuss what will happen if the ERADE is used in such cases. If the urn allocation is the target allocation, and the ERADE (with $\alpha = 1/2$, $m_0 = 2$) is used for treatment allocation, then there will be about 121 patients in the ECMO and 64 patients in the conventional treatment, on average. Based on 10,000 simulated trials, the variability $\sigma^2$ is 0.28, and the expected number of deaths is 74 death, as compared to 82 in the actual trial. The adaptive design utilizes the better treatment more often to save lives.

Next, we compare the variability and power of the ERADE and the RPW rule under the settings of $P_1 = 65/93$ and $P_2 = 38/92$. The expected power under both designs is 0.969. Because they are randomized procedures, the total number of patients assigned to the ECMO group is a random variable. Based on the 10,000 simulated trials, we noticed that in 99 percent of the trials under the ERADE there were more than 52 patients assigned to the conventional treatment group, for a power of 0.941 or higher. Under the RPW rule, only 39 or fewer patients were assigned to the conventional treatment in 1 percent of the trials, for a power of 0.904 or less. Also, based on the 10,000 simulated trials, the ERADE always assigns more patients to the ECMO group; however, the RPW rule assigned more patients to the conventional group in 114 trials. The number of patients assigned to the ECMO group for the ERADE and the RPW are summarized in the box plots of Figure 1. Even at the sample size 185, we can see the advantage of using the proposed ERADE over the randomized player–the-winner rule.

In practice, one may modify the RPW by following changes: (i) first assign $m_0$ (say $m_0 = 2$) subjects to each treatment by using restricted randomization; (ii) update the urn based on their (first $2m_0$ subjects) responses according the RPW rule (add one same type ball if a success, add one opposite type ball if a failure); (iii) use the RPW rule after the first $2m_0$ subjects. For a small size clinical trial (e.g., the Michigan ECMO study), the proposed modification will reduce the variability of the RPW rule. For relatively large size trial (e.g., the UK ECMO trial), the modification with a small $m_0$ (say $m_0 = 2$) will not reduce the variability significantly. However, the modification with a large $m_0$ is not recommended, because it will push the RPW rule to balance the allocation proportion. Also, this modification does not change the asymptotic properties of RPW. It is worthwhile to point out that urn models usually have much larger asymptotic variability than their corresponding ERADEs (with the same limiting allocation proportion) by comparing Theorem 2.1 with the results of Bai and Hu (2005) and Zhang, Hu and Cheung (2006).

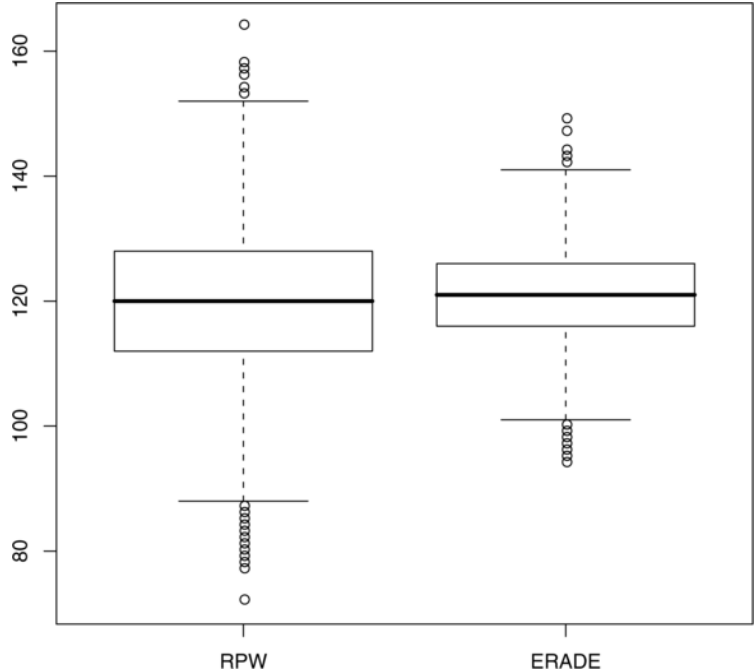

Fig. 1. *Box plots for the number of patients assigned to the ECMO group based on 10,000 simulated trials.*

**5. Conclusion remarks.** We have proposed a family of response-adaptive designs that are fully randomized and asymptotically efficient. The ERADE can be viewed as a generalization of the Efron's biased coin design for any desired allocation function, which may depend on the unknown parameters. The asymptotic properties derived here provide the theoretical foundation for inferences based on the ERADE. They also help to understand why Efron's biased coin design is efficient for balanced trials. The examples in Sections 3 and 4 provide a guideline for the use of the proposed ERADE in practice.

In this paper, we have assumed that the responses in each treatment group are available without delay. In practice, there is no logistical difficulty in incorporating delayed responses into the ERADE, provided that some responses become available during the course of the allocation in the experiment; thus, we can always update the estimates whenever new data become available. For clinical trials with uniform (or exponential) patient entry and exponential response times [see, e.g., Bai, Hu and Rosenberger (2002)], it is easy to verify that both Theorems 2.1 and 2.2 remain true.

In Theorem 2.1, we assumed that the parameter estimate has Bahadur-type representation as in (2.1). This representation is usually satisfied in applications (e.g., moment estimation, maximum likelihood estimation, estimate based on estimating equation, etc.). Theorem 2.1 ensures the consistency and asymptotical normality of the parameter estimate after using ERADE in clinical trial. Baldi Antognini and Giovagnoli (2005) have considered the properties of maximum likelihood estimation for some other response-adaptive designs. Also, we have focused on ERADE for two treatments. Recently, Tymofyeyev, Rosenberger and Hu (2007) proposed a general framework to obtain optimal allocation proportions for comparing two or more treatments. How to generalize ERADE to $K$-treatment case is a new research topic.

In some clinical trials, it is important to take covariates into consideration. To balance known covariates in two treatment groups, Pocock and Simon (1975) proposed a covariate-adaptive randomization procedure by using Efron's biased coin design. Since an important goal of clinical trials is to estimate a treatment effect, one may assign subjects to treatments according to their covariates for efficiency [Atkinson (1982)]. In some clinical trials, one may assign subjects to treatments according to their covariates for ethical reasons [see, e.g., Zhang et al. (2007)]. The ERADE could play an important role in this area; further research beyond the work in the present paper is clearly needed to account the effect of covariates in response-adaptive designs.

## APPENDIX: PROOFS

To prove the main theorems, we first prove the following lemma.

LEMMA A.1. *Suppose* $\{M_n, \mathcal{F}_n; n \geq 1\}$ *is a martingale with* $\mathsf{E}|\Delta M_n|^2 \leq C_0$, *where* $\Delta M_n = M_n - M_{n-1}$ *and* $\{\mathcal{F}_n\}$ *is a filter of sigma fields with* $\mathcal{F}_1 \subset \mathcal{F}_2 \subset \cdots$. *Then, for any* $L > 4$,

$$(A.1) \qquad E\left[\max_{m \leq L} |M_n - M_{n-m}|\right] \leq \sqrt{C_0 L},$$

$$(A.2) \qquad E\left[\max_{m: L \leq m \leq n} \frac{|M_n - M_{n-m}|}{m}\right] \leq 3\sqrt{C_0/L}.$$

PROOF. It is easy to show (A.1). For (A.2), we have

$$E\left[\max_{m:L\leq m\leq n}\frac{|M_n - M_{n-m}|}{m}\right]$$
$$\leq \sum_{j:\log_2 L\leq j\leq \log_2 n} E\left[\max_{2^{j-1}\leq m\leq 2^j}\frac{|M_n - M_{n-m}|}{2^{j-1}}\right]$$
$$\leq \sum_{j:j\geq \log_2 L} \frac{\sqrt{C_0 2^j}}{2^{j-1}} \leq 3\sqrt{C_0/L}.$$

$\square$

PROOF OF THEOREM 2.1. Let $\mathcal{F}_m = \sigma(X_{1,1},\ldots,X_{m,1},\boldsymbol{\xi}_1,\ldots,\boldsymbol{\xi}_m)$ be the sigma field generated by previous $m$ stages. Then, under $\mathcal{F}_{m-1}$, $X_{m,1}$ and $\boldsymbol{\xi}_m$ are independent, and $\mathsf{E}[X_{m,1}|\mathcal{F}_{m-1}] = p_m$. Let $M_n = \sum_{m=1}^n \Delta M_m$, where $\Delta M_m = X_{m,1} - \mathsf{E}[X_{m,1}|\mathcal{F}_{m-1}]$. Then, $\{\Delta M_m; m\geq 1\}$ is a sequence of bounded martingale differences with $|\Delta M_m| \leq 1$, and

$$N_{n,1} = N_{n-1,1} + p_n + \Delta M_n.$$

Write

$$U_n = \sum_{m=0}^{n-1} \alpha\widehat{\rho}_m + M_n - n\widehat{\rho}_n.$$

Then,

(A.3) $$N_{n,1} - n\widehat{\rho}_n \leq N_{n-1,1} - (n-1)\widehat{\rho}_{n-1} + \Delta U_n$$
$$\text{if } N_{n-1,1} - (n-1)\widehat{\rho}_{n-1} > 0,$$

where $\Delta U_n = U_n - U_{n-1}$. Let $l = l_n = \max\{m: 2m_0 + 1 \leq m \leq n, N_{m,1} - m\widehat{\rho}_m \leq 0\}$, where $\max \varnothing = 2m_0$. According to (A.3),

$$N_{n,1} - n\widehat{\rho}_n \leq \Delta U_n + \Delta U_{n-1} + \cdots + \Delta U_{l+1} + N_{l,1} - l\widehat{\rho}_l$$
$$= U_n - U_l + N_{l,1} - l\widehat{\rho}_l.$$

Notice that $N_{l,1} - l\widehat{\rho}_l \leq 0$ if $l \geq 2m_0 + 1$, and $N_{l,1} - l\widehat{\rho}_l \leq 2m_0$ if $l = 2m_0$. So, in either case

(A.4) $$N_{n,1} - n\widehat{\rho}_n \leq U_n - U_l + 2m_0.$$

We will show that

(A.5) $U_n - U_{l_n} \leq o_P(\sqrt{n})$ and $U_n - U_{l_n} \leq O(\sqrt{n\log\log n})$ a.s.

If (A.5) is proved, then

$$N_{n,1} - n\widehat{\rho}_n \leq o_P(\sqrt{n}) \quad \text{and} \quad N_{n,1} - n\widehat{\rho}_n \leq O(\sqrt{n\log\log n}) \quad \text{a.s.,}$$

and, by symmetry,
$$(n - N_{n,1}) - n(1 - \widehat{\rho}_n) \leq o_P(\sqrt{n}),$$
$$(n - N_{n,1}) - n(1 - \widehat{\rho}_n) \leq O(\sqrt{n \log \log n}) \quad \text{a.s.}$$

(2.11) is proved according to (A.4) and (A.5).

For proving (A.5), we first show that $N_{n,1}/n \to v$ a.s. If $N_{n,1} \to \infty$ (resp. $N_{n,2} = n - N_{n,1} \to \infty$), then $\widehat{\boldsymbol{\theta}}_{n,1} \to \boldsymbol{\theta}_1$ (resp. $\widehat{\boldsymbol{\theta}}_{n,2} \to \boldsymbol{\theta}_2$) by Lemma A.4 of Hu and Zhang (2004). If $\sup_n N_{n,1} < \infty$ [resp. $\sup_n(n - N_{n,1}) < \infty$], then the value of $\widehat{\boldsymbol{\theta}}_{n,1}$ (resp. $\widehat{\boldsymbol{\theta}}_{n,2}$) will not change when $n$ is large enough. In either case, $\widehat{\boldsymbol{\theta}}_{n,1}$ and $\widehat{\boldsymbol{\theta}}_{n,2}$ have finite limits. It follows that there is a $u$ with $0 < u < 1$, such that

$$(A.6) \qquad \widehat{\rho}_n \to u \quad \text{a.s.}$$

by the continuity of $\rho(\cdot)$. On the other hand, according to the law of large numbers for martingales,

$$(A.7) \qquad M_n = o(n) \quad \text{a.s.}$$

Combining (A.6) and (A.7) yields $U_n \sim -n(1-\alpha)u$ a.s. And then, by (A.4), $(N_{n,1}/n - \widehat{\rho}_n)^+ \to 0$ a.s. Similarly, $((1 - N_{n,1}/n) - (1 - \widehat{\rho}_n))^+ \to 0$ a.s. It follows that $N_{n,1}/n - \widehat{\rho}_n \to 0$ a.s., which, together with (A.6), implies

$$(A.8) \qquad \lim_{n \to \infty} \frac{N_{n,1}}{n} = \lim_{n \to \infty} \widehat{\rho}_n = u \quad \text{a.s.}$$

On the other hand, (A.8) implies that $N_{n,1} \to \infty$ and $n - N_{n,1} \to \infty$. By Lemma A.4 of Hu and Zhang (2004), again, we have $\widehat{\boldsymbol{\Theta}}_n \to \boldsymbol{\Theta}$ a.s. and

$$(A.9) \qquad \widehat{\boldsymbol{\Theta}}_n - \boldsymbol{\Theta} = O\left(\sqrt{\frac{\log \log n}{n}}\right) \quad \text{a.s.}$$

So, the limit $u$ in (A.6) and (A.8) must be $v = \rho(\boldsymbol{\Theta})$, according to the continuity of $\rho(\cdot)$.

Now,
$$U_n - U_l = -(n-l)(1-\alpha)v + \sum_{m=l}^{n-1} \alpha(\widehat{\rho}_m - v) - (n-l)(\widehat{\rho}_n - v)$$
$$(A.10) \qquad + M_n - M_l + l(\widehat{\rho}_l - \widehat{\rho}_n)$$
$$= (n-l)[-(1-\alpha)v + o(1)] + M_n - M_l + l(\widehat{\rho}_l - \widehat{\rho}_n) \quad \text{a.s.}$$

Also, according to (2.1) and the fact that $N_{n,1}/n \to v$ a.s., it is easily shown that $\sqrt{n}(\widehat{\boldsymbol{\Theta}}_n - \boldsymbol{\Theta}) = O_P(1)$. So, we have

$$(A.11) \quad \widehat{\rho}_n - v = \rho(\widehat{\boldsymbol{\Theta}}_n) - \rho(\boldsymbol{\Theta}) = (\widehat{\boldsymbol{\Theta}}_n - \boldsymbol{\Theta})\frac{\partial \rho}{\partial \mathbf{y}}\bigg|_{\boldsymbol{\Theta}} + O(\|\widehat{\boldsymbol{\Theta}}_n - \boldsymbol{\Theta}\|^2)$$

$$(A.12) \qquad = (\widehat{\boldsymbol{\Theta}}_n - \boldsymbol{\Theta})\frac{\partial \rho}{\partial \mathbf{y}}\bigg|_{\boldsymbol{\Theta}} + O_P\left(\frac{1}{n}\right).$$

From (2.1), it follows that

$$|l(\widehat{\rho}_l - \widehat{\rho}_n)| \leq Cl\|\widehat{\boldsymbol{\Theta}}_n - \widehat{\boldsymbol{\Theta}}_l\| + O_P(1)$$

$$\leq C \sum_{k=1}^{2} \frac{l}{N_{n,k}} \left\| \sum_{m=l+1}^{n} X_{m,k}(\boldsymbol{\xi}_{m,k} - \boldsymbol{\theta}_k) \right\|$$

$$+ C \sum_{k=1}^{2} \frac{l|N_{l,k} - N_{n,k}|}{N_{n,k}} \left\| \frac{1}{N_{l,k}} \sum_{m=1}^{l} X_{m,k}(\boldsymbol{\xi}_{m,k} - \boldsymbol{\theta}_k) \right\| + o_P(\sqrt{n})$$

$$= o(1) \cdot (n-l) + C\|\mathbf{Q}_n - \mathbf{Q}_l\| + o_P(\sqrt{n}),$$

where $\mathbf{Q}_{n,k} = \sum_{m=1}^{n} X_{m,k}(\boldsymbol{\xi}_{m,k} - \boldsymbol{\theta}_k)$ and $\mathbf{Q}_n = (\mathbf{Q}_{n,1}, \mathbf{Q}_{n,2})$, $k = 1, 2$. Notice that $\{\mathbf{Q}_n\}$ and $\{M_n\}$ are both martingales. Notice that, for any $1 \leq L \leq n$, we have

$$\|\mathbf{Q}_n - \mathbf{Q}_l\| \leq (n-l) \cdot \max_{L \leq m \leq n} \frac{\|\mathbf{Q}_n - \mathbf{Q}_{n-m}\|}{m} + \max_{m \leq L} \|\mathbf{Q}_n - \mathbf{Q}_{n-m}\|.$$

Now, choose $L = L(n) \to \infty$ with $L = o(n)$. Then, by Lemma A.1,

$$\max_{L \leq m \leq n} \frac{\|\mathbf{Q}_n - \mathbf{Q}_{n-m}\|}{m} = O_P(\sqrt{1/L}) = o_P(1)$$

and

$$\max_{m \leq L} \|\mathbf{Q}_n - \mathbf{Q}_{n-m}\| = O_P(\sqrt{L}) \leq o_P(\sqrt{n}).$$

It follows that

$$\|\mathbf{Q}_n - \mathbf{Q}_l\| = o_P(1) \cdot (n-l) + o_P(\sqrt{n}).$$

Similarly,

$$|M_n - M_l| = o_P(1) \cdot (n-l) + o_P(\sqrt{n}).$$

We conclude that

$$U_n - U_l \leq (n-l)[-(1-\alpha)v + o_P(1)] + |M_n - M_l| + C\|\mathbf{Q}_n - \mathbf{Q}_l\| + o_P(\sqrt{n})$$

$$\leq (n-l)[-(1-\alpha)v + o_P(1)] + o_P(\sqrt{n}) = o_P(\sqrt{n}).$$

The first part of (A.5) is proved. For the second part, notice that

$$M_n = O(\sqrt{n \log \log n}) \qquad \text{a.s.}$$

due to the law of iterated logarithm for martingales, and

(A.13) $$n(\widehat{\rho}_n - v) = O(\sqrt{n \log \log n}) \qquad \text{a.s.}$$

due to (A.9) and (A.11). The second part of (A.5) is proved by (A.10).

Finally, from Lemma 1 of Hu, Rosenberger and Zhang (2006), we have

(A.14) $$\sqrt{n}(\widehat{\boldsymbol{\Theta}}_n - \boldsymbol{\Theta}) \xrightarrow{D} \text{N}(\mathbf{0}, \mathbf{V}).$$

Together with (2.11) and (A.12), (A.14) yields (2.13). Combining (2.11) and (A.13) yields (2.14). □

**Acknowledgments.** Special thanks go to the anonymous referees, the associate editor and the editors for their comments, which led to a much improved version of this paper.

## REFERENCES

Atkinson, A. C. (1982). Optimal biased coin designs for sequential clinical trials with prognostic factors. *Biometrika* **69** 61–67. MR0655670

Bai, Z. D. and Hu, F. (2005). Asymptotics in randomized urn models. *Ann. Appl. Probab.* **15** 914–940. MR2114994

Bai, Z. D., Hu, F. and Rosenberger, W. F. (2002). Asymptotic properties of adaptive designs with delayed response. *Ann. Statist.* **30** 122–139. MR1892658

Baldi Antognini, A. and Giovagnoli, A. (2004). A new 'biased coin design' for sequential allocation of two treatments. *J. Roy. Statist. Soc. Ser. C* **53** 651–664. MR2087777

Baldi Antognini, A. and Giovagnoli, A. (2005). On the large sample optimality of sequential designs for comparing two or more treatments. *Sequential Anal.* **24** 205–217. MR2154938

Bartlett, R. H., Roloff, D. W., Cornell, R. G., Andrews, A. F., Dillon, P. W. and Zwischenberger, J. B. (1985). Extracorporeal circulation in neonatal respiratory failure: A prospective randomized study. *Pediatrics* **76** 479–487.

Burman, C. F. (1996). On sequential treatment allocation in clinical trials. Dept. Mathematics, Goteborg.

Efron, B. (1971). Forcing a sequential experiment to be balanced. *Biometrika* **58** 403–417. MR0312660

Gwise, T. E., Hu, J. and Hu, F. (2008). Optimal biased coins for two arm clinical trials. *Stat. Interface* **1** 125–136. MR2425350

Hu, F. and Rosenberger, W. F. (2003). Optimality, variability, power: Evaluating response-adaptive randomization procedures for treatment comparisons. *J. Amer. Statist. Assoc.* **98** 671–678. MR2011680

Hu, F. and Rosenberger, W. F. (2006). *The Theory of Response-Adaptive Randomization in Clinical Trials*. Wiley, New York. MR2245329

Hu, F., Rosenberger, W. F. and Zhang, L.-X. (2006). Asymptotically best response-adaptive randomization procedures. *J. Statist. Plann. Inference* **136** 1911–1922. MR2255603

Hu, F. and Zhang, L.-X. (2004). Asymptotic properties of doubly adaptive biased coin designs for multitreatment clinical trials. *Ann. Statist.* **32** 268–301. MR2051008

Ivanova, A. (2003). A play–the-winner type urn model with reduced variability. *Metrika* **58** 1–13. MR1999248

Jennison, C. and Turnbull, B. W. (2000). *Group Sequential Methods with Applications to Clinical Trials*. Chapman & Hall/CRC Press, Boca Raton, FL. MR1710781

Pocock, S. J. and Simon, R. (1975). Sequential treatment assignment with balancing for prognostic factors in the controlled clinical trial. *Biometrics* **31** 103–115.

Rosenberger, W. F. and Lachin, J. M. (2002). *Randomization in Clinical Trials Theory and Practice*. Wiley, New York. MR1914364

Rosenberger, W. F., Stallard, N., Ivanova, A., Harper, C. and Ricks, M. (2001). Optimal adaptive designs for binary response trials. *Biometrics* **57** 909–913. MR1863454

Smith, R. L. (1984). Properties of Biased coin designs in sequential clinical trials. *Ann. Statist.* **12** 1018–1034. MR0751289

Tymofyeyev, Y., Rosenberger, W. F. and Hu, F. (2007). Implementing optimal allocation in sequential binary response experiments. *J. Amer. Statist. Assoc.* **102** 224–234. MR2345540

UK COLLABORATIVE ECMO TRIAL GROUP (1996). UK collaborative randomised trial of neonatal extracorporeal membrane oxygenation. *Lancet* **348** 75–82.

WARE, J. H. (1989). Investigating therapies of potentially great benefit: ECMO (with discussions). *Statist. Sci.* **4** 298–340. MR1041761

WEI, L. J. (1978). The adaptive biased coin design for sequential experiments. *Ann. Statist.* **6** 92–100. MR0471205

WEI, L. J. and DURHAM, S. (1978). The randomized play-the-winner rule in medical trials. *J. Amer. Statist. Assoc.* **73** 840–843.

ZELEN, M. (1969). Play the winner and the controlled clinical trial. *J. Amer. Statist. Assoc.* **64** 131–146. MR0240938

ZHANG, L. and ROSENBERGER, W. F. (2006). Response-adaptive randomization for clinical trials with continuous outcomes. *Biometrics* **62** 562–569. MR2236838

ZHANG, L.-X., HU, F. and CHEUNG, S. H. (2006). Asymptotic theorems of sequential estimation-adjusted urn models. *Ann. Appl. Probab.* **16** 340–369. MR2209345

ZHANG, L.-X., HU, F., CHEUNG, S. H. and CHAN, W. S. (2007). Asymptotic properties of covariate-adjusted adaptive designs. *Ann. Statist.* **35** 1166–1182. MR2341702

F. HU
DEPARTMENT OF STATISTICS
UNIVERSITY OF VIRGINA
HALSEY HALL
CHARLOTTESVILLE, VIRGINIA 22904-4135
USA
E-MAIL: fh6e@virginia.edu

L.-X. ZHANG
DEPARTMENT OF MATHEMATICS
ZHEJIANG UNIVERSITY
HANGZHOU 310027
PEOPLE'S REPUBLIC OF CHINA
E-MAIL: stazlx@zju.edu.cn

X. HE
DEPARTMENT OF STATISTICS
UNIVERSITY OF ILLINOIS AT
  URBANA–CHAMPAIGN
725 S. WRIGHT
CHAMPAIGN, ILLINOIS 61820
USA
E-MAIL: x-he@uiuc.edu